\RequirePackage{fix-cm}
\documentclass[smallextended]{svjour3}       
\smartqed  
\usepackage{graphicx}
\usepackage{amsmath}
\usepackage{ dsfont }
\usepackage{algorithm}
\usepackage{algorithmic}
\usepackage{xcolor,colortbl}
\usepackage{subfigure}
\usepackage{booktabs}
\usepackage{tabularx}
\usepackage{placeins}
\allowdisplaybreaks
\begin{document}

\title{Data Fitting with Signomial Programming Compatible Difference of Convex Functions
}

\titlerunning{Data Fitting with Signomial Programming Compatible DC Functions}        

\author{Cody Karcher
}


\institute{C. Karcher \at
              Department of Aeronautics and Astronautics, Massachusetts Institute of Technology \\ 
              77 Massachusetts Avenue, Cambridge, MA 02139, USA \\
              \email{ckarcher@mit.edu}  \\ 
}

\date{Received: date / Accepted: date}

\maketitle
\begin{abstract}
Signomial Programming (SP) has proven to be a powerful tool for engineering design optimization, striking a balance between the computational efficiency of Geometric Programming (GP) and the extensibility of more general optimization methods like Sequential Quadratic Programming (SQP).  But when an existing engineering analysis tool is incompatible with the mathematics of the SP formulation, options are limited.  Previous literature has suggested schemes for fitting GP compatible models to pre-computed data, but no methods have yet been proposed that take advantage of the increased modeling flexibility available in SP.  This paper describes a new Soft Difference of Max Affine (SDMA) function class that is constructed from existing methods of GP compatible fitting and the theory of Difference of Convex (DC) functions.  When a SDMA function is fit to data in log-log transformed space, it becomes either a signomial or a set of signomials upon inverse transformation.  Three examples of fitting are presented here, including simple test cases in 2D and 3D, and a fit to the performance data of the NACA 24xx family of airfoils.  In each case, RMS error is driven to less than 1\%. 
%
\end{abstract}
\section{Motivation}
\label{motivation}
Geometric Programming (GP) and Signomial Programming (SP) are two related classes of non-linear optimization formulations that have recently been applied, with great success, to the field of aircraft design \cite{hoburg2014geometric,torenbeek2013advanced,hoburg2013fast,kirschen2016signomial,brown2018vehicle,york2018efficient,burton2018solar,lin2020simultaneous,kirschen2018application,york2018turbofan,saab2018robust,hall2018assessment}.  Despite these successes, the GP and SP formulations do not allow the use of black box analysis tools that are common in practical design problems \cite{martins2013multidisciplinary}.  In aircraft design, these tools include high fidelity analyses such as Computational Fluid Dynamics (CFD) and Finite Element Analysis (FEA), and low fidelity analyses like XFOIL \cite{drela1989xfoil} and AVL \cite{drela2010avl}.  Previous work by Hoburg et. al. \cite{hoburg2016data} addressed this issue for Geometric Programming by proposing three methods for fitting a GP compatible function to a pre-computed dataset, but the additional modeling flexibility available in Signomial Programming goes untreated.  This work fills this gap by proposing a method for fitting a SP compatible function to data in cases where the Hoburg methods \cite{hoburg2016data} lack sufficient accuracy.

\section{Geometric Programming and Signomial Programming}
\label{intro}
\subsection{Geometric Programming}
Geometric Programs are built upon two fundamental building blocks: monomial and posynomial functions.  A monomial function is defined as the product of a leading constant with each variable raised to a real power \cite{boyd2007tutorial}:
\begin{equation}
    m(\textbf{x}) = c{x_1}^{a_1}{x_2}^{a_2}...\;{x_n}^{a_n} = c \prod_{i=1}^{N} x_i^{a_i}
\end{equation}
A posynomial is simply the sum of monomials \cite{boyd2007tutorial}, which can be defined in notation as:
\begin{equation}
  p(\textbf{x}) = m_1(\textbf{x}) +  m_2(\textbf{x}) + ... + m_n(\textbf{x}) = \sum_{k=1}^{K} c_k \prod_{i=1}^{N} x_i^{a_{ik}}
\end{equation}
From these two building blocks, it is possible to construct the definition of a GP in standard form \cite{boyd2007tutorial}:
\begin{equation}
    \begin{aligned}
        & \underset{\mathbf{x}}{\text{minimize}}
        & & p_0 (\textbf{x}) \\
        & \text{subject to}
        & & m_i (\textbf{x}) = 1, \; i = 1, \ldots, N\\
        &&& p_j (\textbf{x}) \leq 1, \; j = 1, \ldots, M
    \end{aligned}
    \label{GP_standard}
\end{equation}
When constraints and objectives can be written in the form specified in Equation \ref{GP_standard} it is said that the problem is \textit{GP compatible}.
\subsection{Signomial Programming} \label{sig_sec}
Signomal Programs (SPs) are a logical extension of Geometric Programs that allow the inclusion of negative leading constants and a broader set of equality constraints.  The key building blocks of the Signomial Programing are signomials, which are the difference between two posynomials $p(\textbf{x})$ and $n(\textbf{x})$:
\begin{equation}
  s(\textbf{x}) = p(\textbf{x}) -  n(\textbf{x}) = \sum_{k=1}^{K} c_k \prod_{i=1}^{N} x_i^{a_{ik}} - \sum_{p=1}^{P} d_p \prod_{i=1}^{N} x_i^{g_{ik}}
\end{equation}
The posynomial $n(\textbf{x})$ is often referred to as a `neginomial' because it is made up of all of the terms with negative leading coefficients.  From this definition, the standard form for a Signomial Program can be constructed \cite{kirschen2016signomial}:
\begin{equation}
    \begin{aligned}
        & \underset{\mathbf{x}}{\text{minimize}}
        & & \frac{p_0 (\textbf{x})}{n_0 (\textbf{x})} \\
        & \text{subject to}
        & & s_i (\textbf{x}) = 0, \; i = 1, \ldots, N\\
        &&& s_j (\textbf{x}) \leq 0, \; j = 1, \ldots, M
    \end{aligned}
    \label{SP_standard}
\end{equation}
however, another useful form is \cite{kirschen2016signomial}:
\begin{equation}
    \begin{aligned}
        & \underset{\mathbf{x}}{\text{minimize}}
        & & \frac{p_0 (\textbf{x})}{n_0 (\textbf{x})} \\
        & \text{subject to}
        & & \frac{p_i (\textbf{x})}{n_i (\textbf{x})} = 1, \; i = 1, \ldots, N\\
        &&& \frac{p_j (\textbf{x})}{n_j (\textbf{x})} \leq 1, \; j = 1, \ldots, M
    \end{aligned}
    \label{SP_standardMod}
\end{equation}
In this alternative form, the neginomial is added to both sides, and then used as a divisor to construct an expression either equal to or constrained by a value of one. 

SPs are \textit{not} convex upon log-log transformation unlike their GP counterparts and therefore must be solved using general non-linear solution methods.  However many signomial programs of interest still benefit from an underlying structure which is well approximated by a log-log convex formulation, and as a result can be efficiently solved by a series of GP approximations via the Difference of Convex Algorithm (DCA).  In this process the various neginomials $n(\textbf{x})$ are replaced with local monomial approximations, yielding substantial benefits over other non-linear solution methods (see \cite{kirschen2018application} and \cite{york2018efficient} for discussion).
\section{Difference of Convex Functions for Data Fitting}
\subsection{Difference of Convex (DC) Functions}
Signomials are extremely versatile functions, far more versatile than posynomials or monomials.  But the construction of signomials is rather unique since they are the difference between two log-log convex functions.

Modern study of Difference of Convex (DC) functions is widely attributed to have started with a paper by Hartman et. al. \cite{hartman1959functions}, which states that a continuous function of bounded variation ($f(\mathbf{x})$) can be written as the difference of two convex functions ($g(\mathbf{x})$ and $h(\mathbf{x})$) \cite{hartman1959functions}:
\begin{equation}
    f(\mathbf{x}) = g(\mathbf{x}) - h(\mathbf{x})
\end{equation}
which clearly mirrors signomial construction.  Further extensions have opened up a wide variety of well known properties that DC functions possess \cite{bavcak2011difference,rudin2012prediction}.  

The implication here is that a significant subset of continuous functions can be decomposed into the difference of two convex functions, including non-differentiable continuous functions.  By extension, it should be possible to fit a DC function to nearly any dataset that can be well approximated with a continuous function.

Hoburg et. al. \cite{hoburg2016data} provided three ways of fitting functions to datasets that were convex under the log-log transformation, but these methods entirely ignore the additional modeling flexibility available in Signomial Programming.  So consider the use of a DC function under the log-log transformation, specifically the difference between two of the functions proposed by Hoburg et. al. \cite{hoburg2016data}.  
\subsection{Function Definitions}
\subsubsection{Notation} 
Consider a data set sampled from a black box mapping from $\mathds{R}^N \rightarrow \mathds{R}$.  Consistent with the notation in Hoburg \cite{hoburg2016data} let the vector $\mathbf{u}_j$ represent the independent variables in $\mathds{R}^N$ for data point $j$ and the scalar $w_j$ represent the output in $\mathds{R}$ for data point $j$.  The logspace variables are then represented as $\mathbf{x}_j = \log \mathbf{u}_j$ and $y_j = \log w_j$.  
\subsubsection{Difference of Max Affine (DMA) Functions} 
The first function proposed by Hoburg et. al. \cite{hoburg2016data} is the Max Affine (MA) function:
\begin{equation}
f_{\text{MA}}(\textbf{x}) = \max_{k=1\ldots K} \left[ b_k + \textbf{a}_k^{\text{T}} \textbf{x} \right]
\label{maxaffine}
\end{equation}

This function class is known to create a convex epigraph.  In fact, any convex function can be reasonably approximated as a max affine function given a sufficient number of affine functions, $K$.  

Now consider the difference between two of these max affine functions (Equation \ref{maxaffine}), which henceforth will be called the Difference of Max Affine (DMA) function:
\begin{equation}
f_{\text{DMA}}(\textbf{x}) = \max_{k=1\ldots K} \left[ b_k + \textbf{a}_k^{\text{T}} \textbf{x} \right] - \max_{m=1\ldots M} \left[ h_m + \textbf{g}_m^{\text{T}} \textbf{x} \right]
\label{differenceofmaxaffine}
\end{equation}

The subtracted term is represented by an entirely separable Max Affine function which is defined by fitting parameters $M$, $h$, and $\mathbf{g}$.  

While the Max Affine function has a realizable transformation back from logspace, the Difference of Max Affine function has no such transformation due to the inability to readily construct a meaningful epigraph or subgraph.  Despite this, the DMA function is quite rapid to fit, and could be of application in other areas where a cheap surrogate is desired for non-convex fitting.  Here, the DMA function is a useful as an intellectual building block, and as a seed for the function class that follows.
\subsubsection{Soft Difference of Max Affine (SDMA) Functions} 
The second function proposed by Hoburg et. al. \cite{hoburg2016data} is the Soft-Max Affine (SMA) function:
\begin{equation}
\label{fsma}
f_{\text{SMA}}(\textbf{x}) = \frac{1}{\alpha} \log \sum_{k=1}^{K} \exp \left( \alpha(b_k + \mathbf{a}_k^{\text{T}} \mathbf{x}) \right)
\end{equation}

The SMA function uses a global softening parameter ($\alpha$) to `smooth' the sharp corners of the Max Affine function and has the benefit of requiring far fewer affine terms $K$ to capture smooth convex functions with reasonable accuracy.  However, the global softening parameter results in a poor representation in regions where the curvature deviates substantially from the global average.

Consider the following function which is the difference between two Soft-Max Affine functions (Equation \ref{fsma}):
\begin{equation}
\label{fsdma}
f_{\text{SDMA}}(\textbf{x}) = \frac{1}{\alpha} \log \sum_{k=1}^{K} \exp \left( \alpha(b_k + \mathbf{a}_k^{\text{T}} \mathbf{x}) \right) - 
                              \frac{1}{\beta}  \log \sum_{m=1}^{M} \exp \left( \beta (h_m + \mathbf{g}_m^{\text{T}} \mathbf{x}) \right)
\end{equation}

which his henceforth referred to as the Soft Difference of Max Affine (SDMA) function.  As with SMA functions when compared to MA functions, the SDMA function compared to the DMA function significantly decreases the number of terms $K+M$ required to approximate smooth functions to a reasonable degree of accuracy by introducing scalar softening parameters $\alpha$ and $\beta$.  

Transforming the SDMA function back to the optimization variables $\mathbf{u}_j$ and $w_j$ proceeds as follows:
\begin{equation}
    \begin{aligned}
        y & = \frac{1}{\alpha} \log \sum_{k=1}^{K} \exp \left( \alpha(b_k + \mathbf{a}_k^{\text{T}} \mathbf{x}) \right) - 
            \frac{1}{\beta}  \log \sum_{m=1}^{M} \exp \left( \beta (h_m + \mathbf{g}_m^{\text{T}} \mathbf{x}) \right) \\
        \exp(y) & = \exp \left( \frac{1}{\alpha} \log \sum_{k=1}^{K} \exp \left( \alpha(b_k + \mathbf{a}_k^{\text{T}} \mathbf{x}) \right) - 
                     \frac{1}{\beta}  \log \sum_{m=1}^{M} \exp \left( \beta (h_m + \mathbf{g}_m^{\text{T}} \mathbf{x}) \right) \right) \\
        w & = \frac{\exp \left( \frac{1}{\alpha} \log \sum_{k=1}^{K} \exp \left( \alpha(b_k + \mathbf{a}_k^{\text{T}} \mathbf{x}) \right) \right)}{\exp \left( \frac{1}{\beta}  \log \sum_{m=1}^{M} \exp \left( \beta (h_m + \mathbf{g}_m^{\text{T}} \mathbf{x}) \right) \right)} \\
        w & = \frac{\exp \left( 
        \log \left( \sum_{k=1}^{K} \exp \left( \alpha(b_k + \mathbf{a}_k^{\text{T}} \mathbf{x}) \right) \right)^{\frac{1}{\alpha}} \right)}{\exp \left( 
        \log \left( \sum_{m=1}^{M} \exp \left( \beta (h_m + \mathbf{g}_m^{\text{T}} \mathbf{x}) \right) \right)^{\frac{1}{\beta} } \right)}
    \end{aligned}
\end{equation}

At this point it is obvious from the definition of a posynomial function that the form will reduce to:
\begin{equation}
w =  \frac{
              \left( \sum_{k=1}^{K} e^{\alpha b_k} \prod_{i=1}^N u_i^{\alpha a_{ik}} \right)^{1/\alpha}   }{
              \left( \sum_{m=1}^{M} e^{\beta  h_m} \prod_{i=1}^N u_i^{\beta g_{im}} \right)^{1/\beta}     }
\label{dsmaReduced}
\end{equation}

This is not compatible with the SP formulation, but consider the following substitutions:
\begin{align}
    p_{convex}  &=  \sum_{k=1}^{K} e^{\alpha b_k} \prod_{i=1}^N u_i^{\alpha a_{ik}} \label{dsmaSet1} \\
    p_{concave} &=  \sum_{m=1}^{M} e^{\beta  h_m} \prod_{i=1}^N u_i^{\beta g_{im}}  \label{dsmaSet2}
\end{align}

Which then reduces Equation \ref{dsmaReduced} to:
\begin{equation}
    w =  \frac{(p_{convex})^{1/\alpha} }{(p_{concave})^{1/\beta} }
    \label{dsmaSet3}
\end{equation}

Thus, taking the three constraints (Equations \ref{dsmaSet1}, \ref{dsmaSet2}, \ref{dsmaSet3}) as a set does result in an SP compatible scheme.  This method of substitution is consistent with other approaches to constructing GP and SP compatible constraints \cite{boyd2007tutorial}.  The table below outlines the operators that should be used in the new constraint set given the original mapping:

\newcolumntype{Y}{>{\centering\arraybackslash}X}
\begin{table}[htbp] 
  \footnotesize 
  \begin{center} 
  \caption{Operators For a Given Initial Constraint} 
  \begin{tabularx}{0.5\textwidth}{  *{4}{Y}} 
      \cmidrule[\heavyrulewidth]{2-4}
       & \multicolumn{3}{c}{Constraint Operator} \\
      \cmidrule{1-4}
      Original & $=$ & $\leq$ & $\geq$ \\
      \cmidrule[\heavyrulewidth]{1-4}
      $w$            &  $=$ & $\leq$ & $\geq$ \\ 
      $p_{convex}$   &  $=$ & $\geq$ & $\leq$ \\ 
      $p_{concave}$  &  $=$ & $\leq$ & $\geq$ \\ \bottomrule
  \end{tabularx} 
 \end{center} 
\end{table}

\subsubsection{Consideration of Implicit Soft Difference of Max Affine (ISDMA) Functions}
Hoburg et. al. \cite{hoburg2016data} proposes a third function class, Implicit Soft Max Affine (ISMA).  Unlike MA and SMA, the ISMA function class is an implicit function of $y$:
\begin{align}
\label{fisma}
\tilde{f}_{\text{ISMA}}(\mathbf{x}, y) &= \log \sum_{k=1}^{K} \alpha_k \exp \left( \alpha_k(b_k + \mathbf{a}_k^{\text{T}}\mathbf{x} - y) \right)
\end{align}

Since the ISMA function class provided a fit with lower RMS error than SMA, it is tempting to write an implicit soft difference of max affine function as:
\begin{equation}
\begin{aligned}
\label{fisdma}
\tilde{f}_{\text{ISDMA}}(\mathbf{x}, y) &= \log \sum_{k=1}^{K} \alpha_k \exp \left( \alpha_k(b_k + \mathbf{a}_k^{\text{T}}\mathbf{x} - y_{convex}) \right) \\
                                      & - \log \sum_{m=1}^{M} \beta_m  \exp \left( \beta_m (h_m + \mathbf{g}_m^{\text{T}}\mathbf{x} - y_{concave}) \right)
\end{aligned}
\end{equation}

Fitting one of these functions to data would require that the data first be decomposed into convex ($y_{convex}$) and concave ($y_{concave}$), which has no obvious solution a priori.  Fortunately, in practice the SDMA function is highly versatile and captures regions with varying curvature for even low fit orders due to the DC construction, somewhat negating the desire for ISDMA functions in the first place.  However, non-differentiable corners and cusps would benefit from this type of function if the issues could be successfully resolved.

\section{How Functions are Fit}
The following process is used for fitting functions to the data points ($\mathbf{u}_j$,$w_j$):
\begin{enumerate}
    \item Apply the log transformation $(\mathbf{x}_j, y_j)=(\log \mathbf{u}_j, \log w_j)$
    \item Fit a function $f$ to the transformed data such that $y_j \approx f(\mathbf{x}_j)$
    \item Relax the equality $y_j = f(\mathbf{x}_j)$ as appropriate to construct the desired constraint
\end{enumerate}

This process is nearly identical to that proposed by Hoburg et. al. \cite{hoburg2016data}, but two differences exist.  First, the work here does not require convex functions in Step 2.  Second, while Hoburg et. al. \cite{hoburg2016data} could readily relax each equality constraint to construct a convex epigraph in Step 3, there is no similar approach when considering functions that are generally not convex and so selecting the appropriate constraint operator ($=$,$\leq$, or $\geq$) will depend on the constraint being considered.

For a set of $m$ data points, the fitting problem can be written as an unconstrained least squares optimization problem:
\begin{equation}
\underset{\gamma}{\text{minimize}} \hspace{1ex} \sum_{i = 1}^m \left(f(\mathbf{x}_i; \mathbf{\gamma}) - y_i \right)^2
\label{leastsquares}
\end{equation}
\noindent
where the fitting parameters are stacked in the vector $\gamma$. 

Based on the work of Hoburg \cite{hoburg2016data}, all of the functions presented here were fit using a Levenberg-Marquardt (LM) algorithm.  The LM algorithm computes a step size at each iteration, but requires a Jacobian to be computed at each step.  Derivatives are therefore presented for each function class with respect to the fitting parameters in the vector $\gamma$ in the two sections that follow. 

For each SDMA function fit, a DMA function is first fit in order to provide an initial guess for the SDMA fitting algorithm.  The ability to quickly achieve a DMA initial guess is critical to the success of SDMA functions, as starting from random positions rarely yields a good result.

In addition, the gradient based nature of the LM algorithm meant that the fitting process required a number of random restarts to be performed from varying initial guesses.  In general 30-100 random restarts were performed for the work presented here, with the specific number varying with the individual fitting problem being considered.

Non-gradient based methods such as a Genetic Algorithm or Particle Swarm Optimization were considered for fitting, but there was no evidence that this approach would be superior to the LM algorithm established in the literature.  A Genetic Algorithm was actually tested, but had worse performance than the LM method even in simple 2D cases.

\subsection{Derivatives for the DMA Function}
\begin{equation}
\frac{\partial f_{\text{DMA}}}{\partial b_i} = \left\{ 
\begin{tabular}{l}
1, \text{if} $i = \text{argmax}_k \hspace{1ex} b_k + \mathbf{a}_k^{\text{T}} \mathbf{x}$ \\
0, otherwise
\end{tabular}
\right.
\end{equation}
\begin{equation}
\hspace{2ex}
\frac{\partial f_{\text{DMA}}}{\partial \mathbf{a}_i} = \left\{ 
\begin{tabular}{l}
$\mathbf{x}^{\text{T}}$, \text{if} $i = \text{argmax}_k \hspace{1ex} b_k + \mathbf{a}_k^{\text{T}} \mathbf{x}$ \\
0, otherwise
\end{tabular}
\right. 
\end{equation}
\begin{equation}
\frac{\partial f_{\text{DMA}}}{\partial h_i} = \left\{ 
\begin{tabular}{l}
-1, \text{if} $i = \text{argmax}_k \hspace{1ex} b_k + \mathbf{a}_k^{\text{T}} \mathbf{x}$ \\
0, otherwise
\end{tabular}
\right.
\end{equation}
\begin{equation}
\hspace{2ex}
\frac{\partial f_{\text{DMA}}}{\partial \mathbf{g}_i} = \left\{ 
\begin{tabular}{l}
$-\mathbf{x}^{\text{T}}$, \text{if} $i = \text{argmax}_k \hspace{1ex} b_k + \mathbf{a}_k^{\text{T}} \mathbf{x}$ \\
0, otherwise
\end{tabular}
\right.
\end{equation}
\subsection{Derivatives for the SDMA Function} 
\begin{align*}
\frac{\partial f_{\text{SDMA}}}{\partial b_i} &= \frac{\exp({\alpha(b_i + \mathbf{a}_i^{\text{T}}\mathbf{x}))}}
                                                      {\sum_{k = 1}^K \exp(\alpha(b_k + \mathbf{a}_k^{\text{T}}\mathbf{x}))} \\[0.06in]
\frac{\partial f_{\text{SDMA}}}{\partial \mathbf{a}_i} &= \frac{\mathbf{x}^{\text{T}} \cdot \exp(\alpha(b_i + \mathbf{a}_i^{\text{T}}\mathbf{x}))}
                                                               {\sum_{k = 1}^K \exp(\alpha(b_k + \mathbf{a}_k^{\text{T}}\mathbf{x}))} \\[0.06in]
\frac{\partial f_{\text{SDMA}}}{\partial \alpha} &= \frac{1}{\alpha} \frac{\sum_{k=1}^K(b_k + \mathbf{a}_k^{\text{T}} \mathbf{x}) \exp(\alpha(b_k + \mathbf{a}_k^{\text{T}}\mathbf{x}))}
                                                                          {\sum_{k = 1}^K \exp(\alpha(b_k + \mathbf{a}_k^{\text{T}}\mathbf{x}))} \\
                                                                         & - \frac{1}{\alpha^2} \log \sum_{k=1}^{K} \exp \left( \alpha(b_k + \mathbf{a}_k^{\text{T}} \mathbf{x}) \right) \\[0.06in]
\frac{\partial f_{\text{SDMA}}}{\partial (\log \alpha)} &= \frac{\sum_{k=1}^K(b_k + \mathbf{a}_k^{\text{T}} \mathbf{x}) \exp(\alpha(b_k + \mathbf{a}_k^{\text{T}}\mathbf{x}))}
                                                                {\sum_{k = 1}^K \exp(\alpha(b_k + \mathbf{a}_k^{\text{T}}\mathbf{x}))} \\
                                                                & - \frac{1}{\alpha} \log \sum_{k=1}^{K} \exp \left( \alpha(b_k + \mathbf{a}_k^{\text{T}} \mathbf{x}) \right) \\[0.06in]
\frac{\partial f_{\text{SDMA}}}{\partial h_i} &= \frac{-\exp({\beta(h_i + \mathbf{g}_i^{\text{T}}\mathbf{x}))}}
                                                      {\sum_{m = 1}^M \exp(\beta(h_m + \mathbf{g}_m^{\text{T}}\mathbf{x}))} \\[0.06in] \stepcounter{equation}\tag{\theequation}\label{sdmaJacobian}
\frac{\partial f_{\text{SDMA}}}{\partial \mathbf{g}_i} &= \frac{-\mathbf{x}^{\text{T}} \cdot \exp(\beta(h_i + \mathbf{g}_i^{\text{T}}\mathbf{x}))}
                                                               {\sum_{m = 1}^M \exp(\beta(h_m + \mathbf{g}_m^{\text{T}}\mathbf{x}))} \\[0.06in]
\frac{\partial f_{\text{SDMA}}}{\partial \beta} &= -\frac{1}{\beta} \frac{\sum_{m=1}^M(h_m + \mathbf{g}_m^{\text{T}} \mathbf{x}) \exp(\beta(h_m + \mathbf{g}_m^{\text{T}}\mathbf{x}))}
                                                                          {\sum_{m = 1}^M \exp(\beta(h_m + \mathbf{g}_m^{\text{T}}\mathbf{x}))} \\
                                                                          & ~~~ + \frac{1}{\beta^2} \log \sum_{m=1}^{M} \exp \left( \beta(h_m + \mathbf{g}_m^{\text{T}} \mathbf{x}) \right) \\[0.06in]
\frac{\partial f_{\text{SDMA}}}{\partial (\log \beta)} &= -\frac{\sum_{m=1}^M(h_m + \mathbf{g}_m^{\text{T}} \mathbf{x}) \exp(\beta(h_m + \mathbf{g}_m^{\text{T}}\mathbf{x}))}
                                                                {\sum_{m = 1}^M \exp(\beta(h_m + \mathbf{g}_m^{\text{T}}\mathbf{x}))} \\
                                                                & ~~~ + \frac{1}{\beta} \log \sum_{m=1}^{M} \exp \left( \beta(h_m + \mathbf{g}_m^{\text{T}} \mathbf{x}) \right)
\end{align*}

\section{Demonstrations of the Fitting Models}
\subsection{A 2D Fitting Problem}
Consider the 2D function, which uses the log transformed variables $\mathbf{x}$ and $y$:
\begin{equation}
    y = \max \left[ -6x - 6, x^4-3x^2 \right]
\end{equation}

This function was selected to test the SDMA function at a non-differentiable point and over variations in curvature.  A set of 101 data points were sampled on even intervals over the function from $x=[-2,2]$, and 5th order MA, SMA, ISMA, and SDMA functions were fit to this data:  
\begin{figure}[htb]
\centering     
\vspace{-0.2in}
\includegraphics[width=0.61\textwidth]{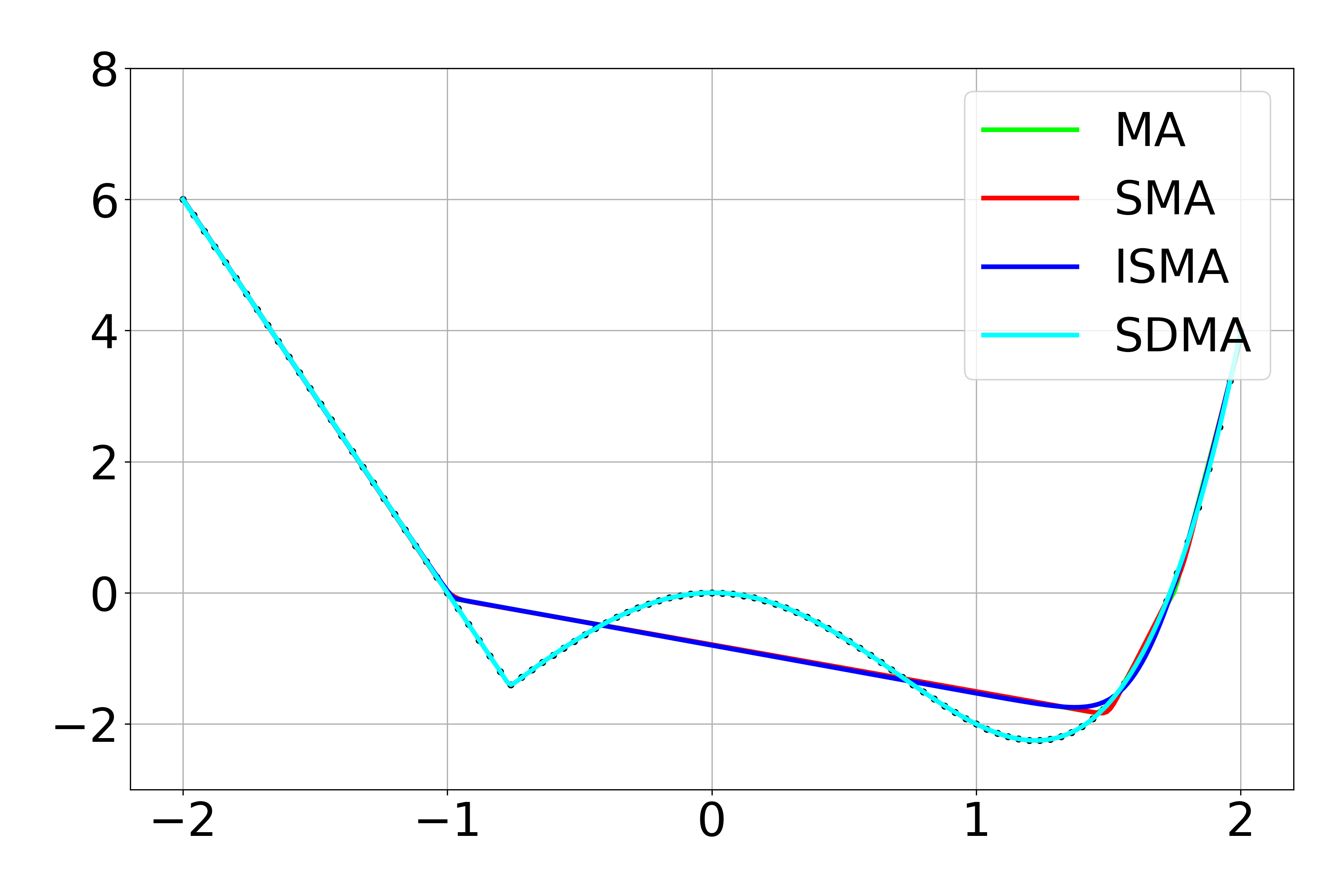}
\vspace{-0.1in}
\caption{Fifth order fits to the simple 2D function.  RMS Error is as follows:  MA: 44.6\%, SMA: 44.6\%, ISMA: 44.8\%, SDMA: 0.149\%}
\label{2dfit}
\end{figure}
\FloatBarrier

\vspace{-0.175in}
The convex fitting methods all converge to nearly identical representations, none of which capture the significant non-convexity to any degree.  In contrast, the SDMA function captures nearly all of the non-convexity in the function, to an RMS error of $1.5\times10^{-3}$.  The key takeaway here is that unlike SMA functions, which can only capture a single curvature due to the single parameter $\alpha$, the SDMA function can capture complex, multi-radius curvature as a direct result of the DC construction.  

Increasing fit order also improves the RMS error as expected: 
\begin{figure}[htb]
\centering     
\includegraphics[width=0.61\textwidth]{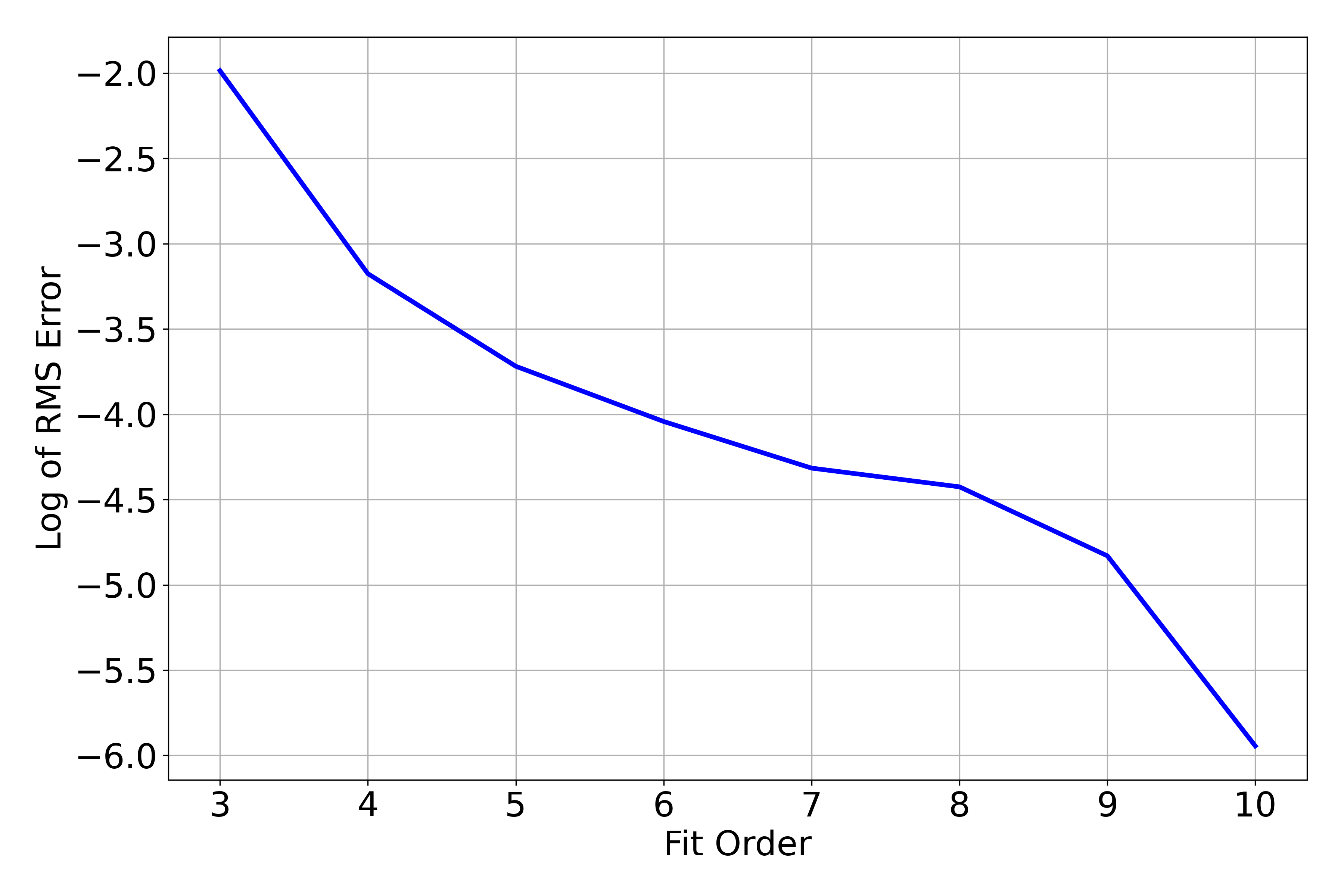}
\vspace{-0.1in}
\caption{Variation of RMS error with increasing fit order}
\label{2dfitError}
\end{figure}
\FloatBarrier

\subsection{A 3D Fitting Problem}
After proving on a 2D case, the next logical step was to consider a 3D dataset.  Similar to the 2D test case, it was desirable to have a non-differentiable region and complex 3D curvature.  Consider this test function, the eigenfunction of the wave equation:
\begin{figure}[htb]
\centering     
\vspace{-0.2in}
\includegraphics[width=0.65\textwidth]{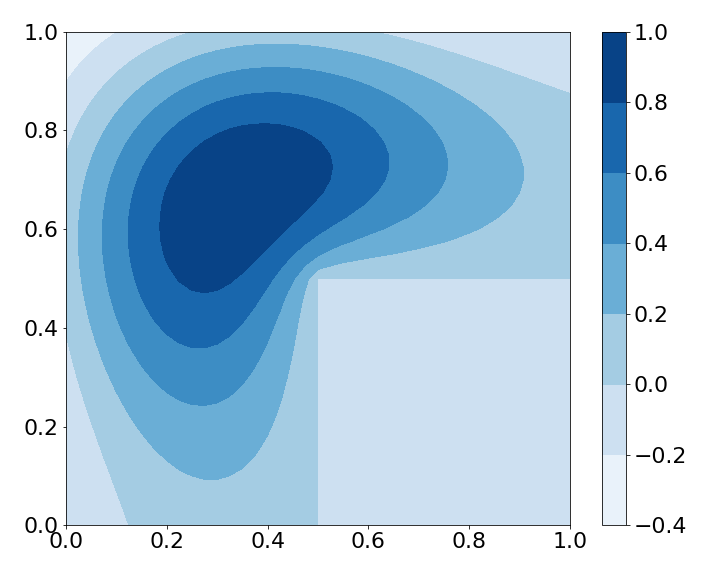}
\vspace{-0.1in}
\caption{The eigenfunction of the wave equation}
\label{matlablogo}
\end{figure}
\FloatBarrier
\vspace{-0.15in}
\noindent
This function has become well known as the Matlab logo.

An SDMA fit yields a reasonably accurate surrogate (consider here 9th order):
\begin{figure}[htb]
\centering     
\subfigure[The fitted 9th order function]{\includegraphics[width=0.47\textwidth]{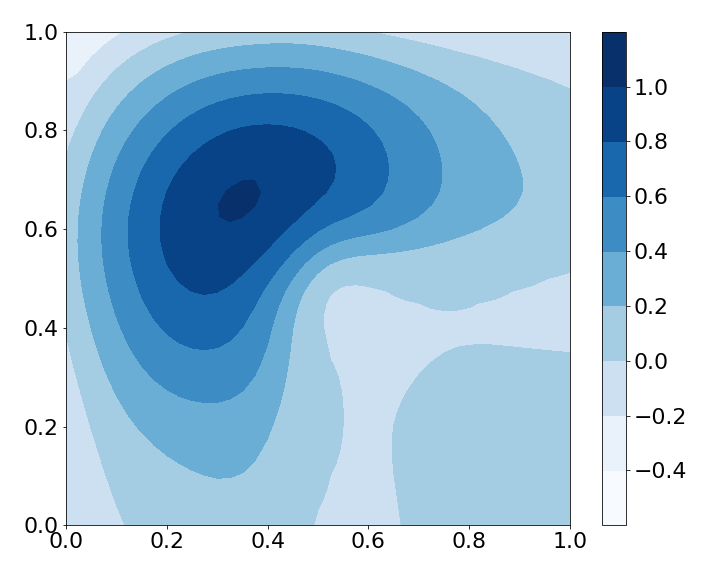}}
\subfigure[Error compared to the true function \label{fig9b}]{\includegraphics[width=0.47\textwidth]{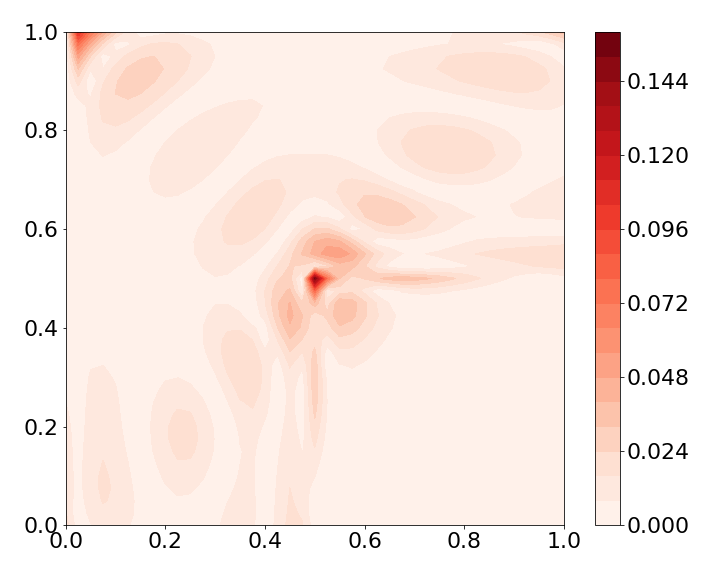}}
\vspace{-0.05in}
\caption{The fitted 9th order function and error.  Note the increased value at the top of the color axis compared to 8th order, indicating a larger maximum error}
\label{matlab9error}
\end{figure}
\FloatBarrier

As might be expected, the fitting scheme struggles along the non-differentiable L-shaped curve, and at the point of the L specifically.  Error near the point (0,1) is due to an inflection in the data where curvature changes from concave to convex in a very small region, and this is difficult to capture in the fit.

And as before, plotting RMS error as a function of increasing fit order shows significant improvement for higher order fits:
\begin{figure}[htb]
\centering     
\includegraphics[width=0.8\textwidth]{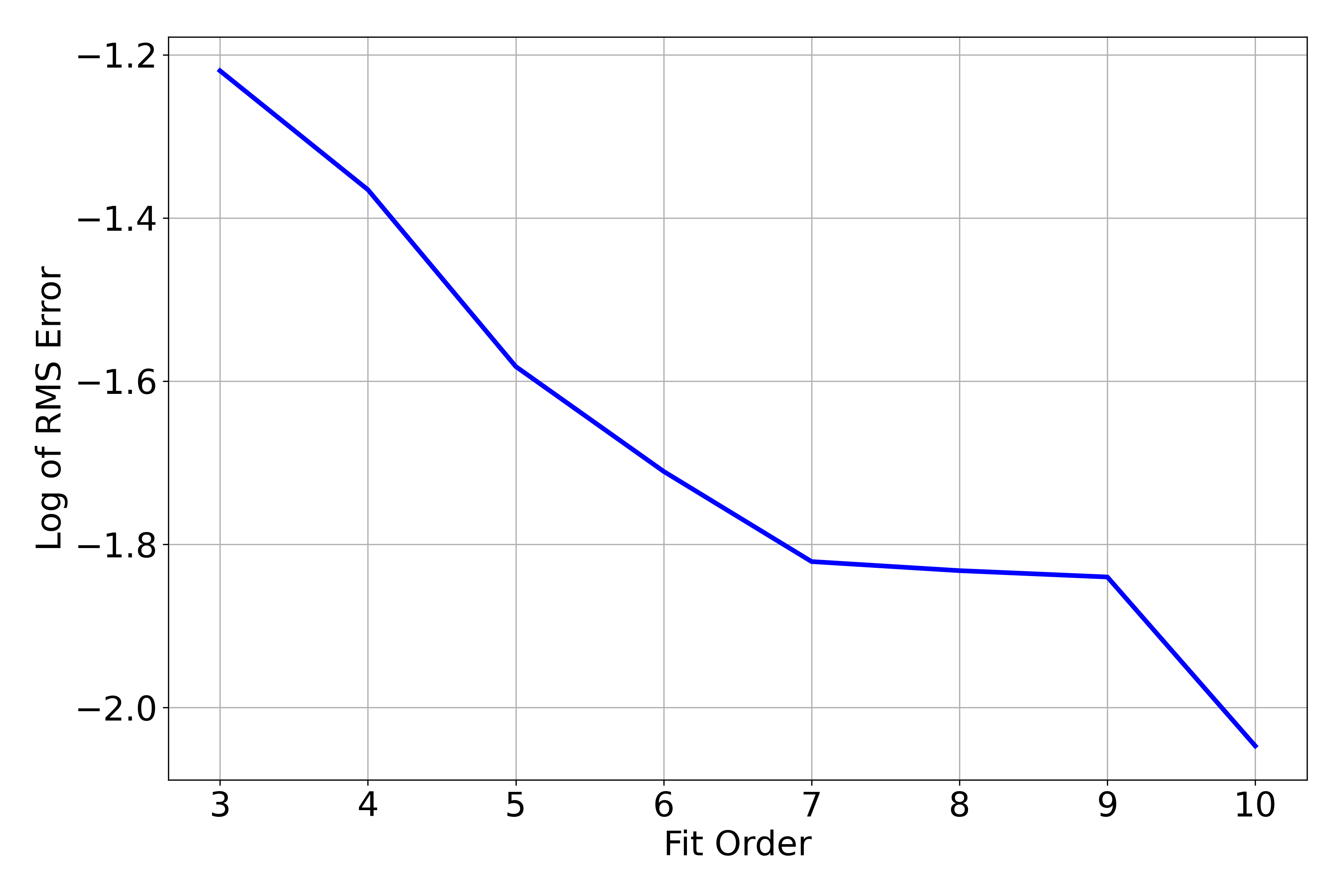}
\caption{Variation of RMS error with increasing fit order}
\label{matlabFitError}
\end{figure}
\FloatBarrier
\subsection{Fitting XFOIL Performance Data of the NACA 24xx Family of Airfoils}
Given the success of the two test cases, the final step was to validate on a more representative engineering example.  The Hoburg et. al. \cite{hoburg2016data} work was able to fit performance data for NACA 24xx airfoils generated from XFOIL \cite{drela1989xfoil}, but only by considering curves of lift coefficient vs. drag coefficient.  In many cases, it is more useful to have two separate curves of lift coefficient vs. angle of attack and drag coefficient vs. angle of attack, but the CL vs. $\alpha$ curve is not compatible with log-log convex fitting techniques.

Consider the problem of fitting the following function:
\begin{equation}
    C_L = f(\alpha,Re,\tau)
\end{equation}

where $C_L$ is the airfoil lift coefficient, $\alpha$ is the airfoil angle of attack, $Re$ is the Reynolds number, and $\tau$ is the airfoil thickness (ex, a $\tau=0.12$ would be a NACA 2412 airfoil).  

XFOIL was used to generate a grid of data from $\alpha=[1,23]$, $Re=[10^6,10^7]$, and $\tau=[.09,.21]$.  This data was then fit with an SDMA function of fit orders varying from 3 to 10.  Plotting below are the results of a 5th order fit:

\begin{figure}[htb]
\centering   
\subfigure[Reynolds Number of $10^6$]{\includegraphics[width=0.49\textwidth]{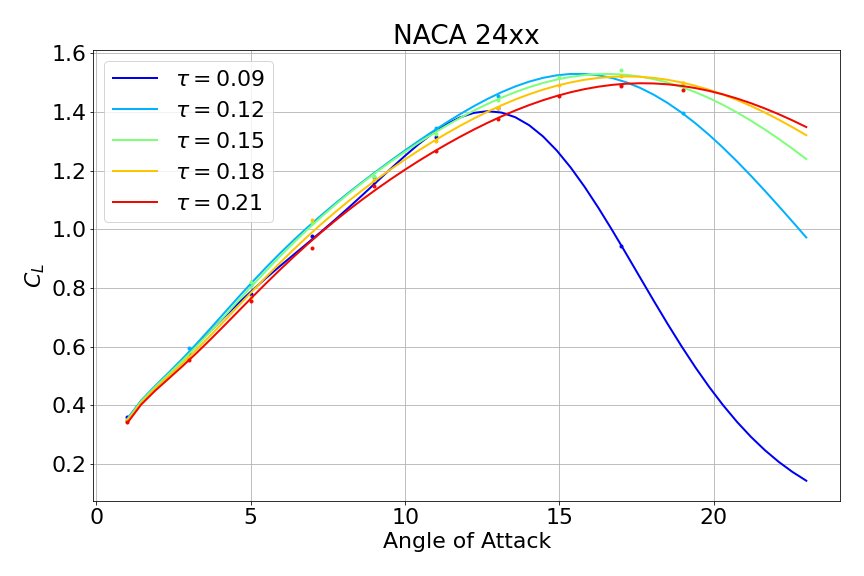}} 
\subfigure[Reynolds Number of $10^{6.5}$]{\includegraphics[width=0.49\textwidth]{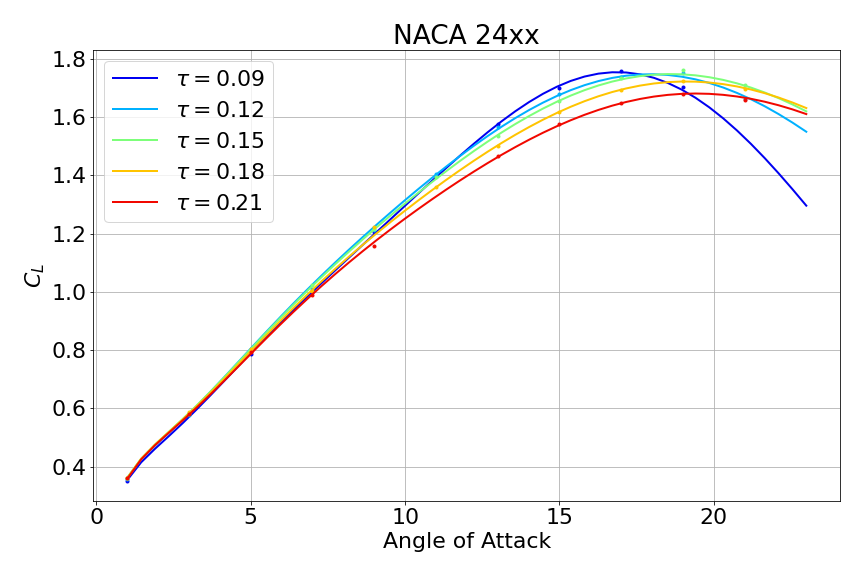}}
\subfigure[Reynolds Number of $10^{7}$]{\includegraphics[width=0.49\textwidth]{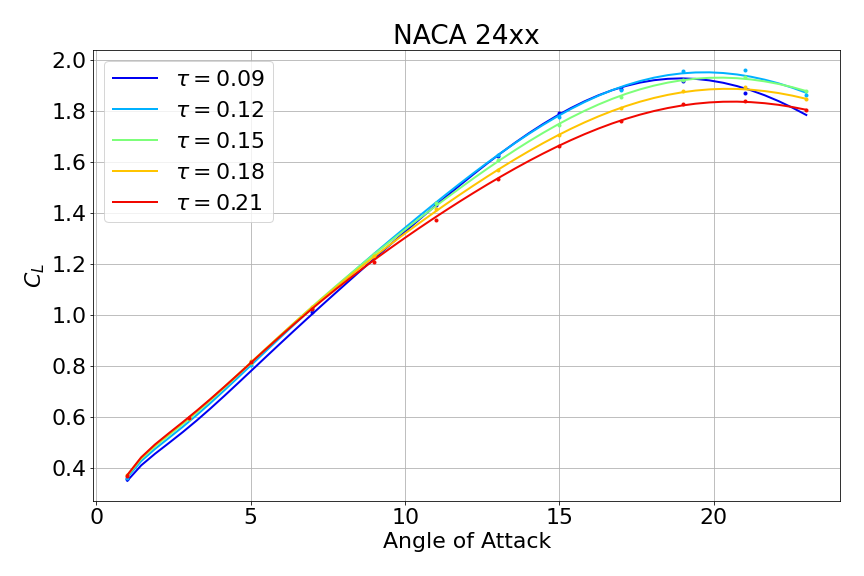}}
\caption{Curves of Lift Coefficient vs. Angle of Attack for the NACA 24xx family of airfoils}
\label{nacaDataFit}
\vspace{-0.15in}
\end{figure}
\FloatBarrier

As before, plotting the RMS error vs. increasing fit order shows improvement with an increased number of fit terms:
\begin{figure}[htb]
\centering     
\includegraphics[width=0.8\textwidth]{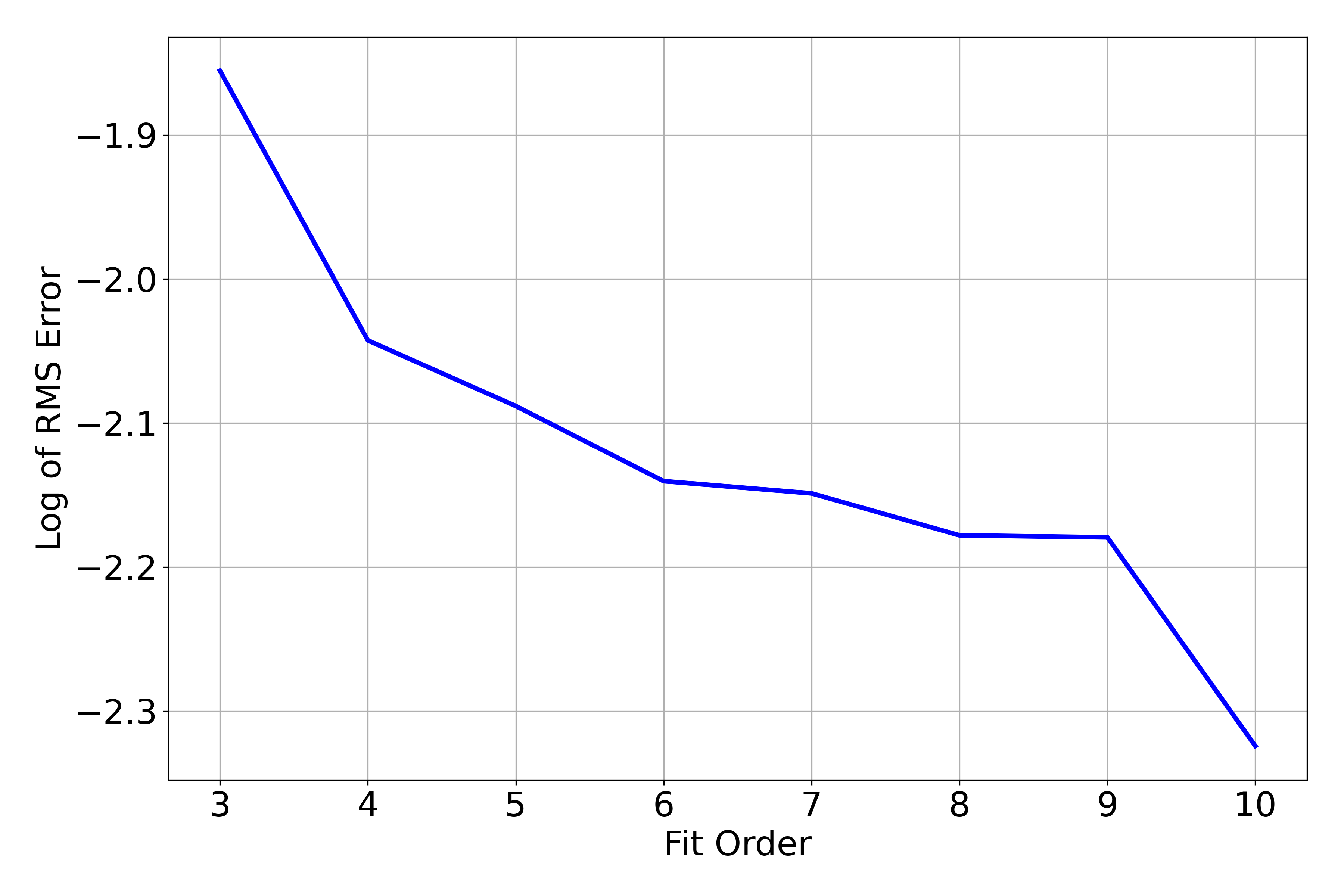}
\vspace{-0.1in}
\caption{Variation of RMS error with increasing fit order}
\label{clFitError}
\end{figure}
\FloatBarrier

\section{Conclusions}
This work serves as a successful development and validation of a Signomial Programming compatible method for fitting black box data.  The Soft Difference of Max Affine (SDMA) function is a crucial link in the chain from Geometric Programming, where low fidelity analysis methods must be used, to Sequential Quadratic Programming or the more recently proposed Logspace Sequential Quadratic Programming \cite{karcher2021logspace}, which typically impose no bounds on the type of analysis that can be utilized.

\begin{acknowledgements}
This material is based on research sponsored by the U.S. Air Force under agreement number FA8650-20-2-2002. The U.S. Government is authorized to reproduce and distribute reprints for Governmental purposes notwithstanding any copyright notation thereon.  The views and conclusions contained herein are those of the authors and should not be interpreted as necessarily representing the official policies or endorsements, either expressed or implied, of the U.S. Air Force or the U.S. Government.

The author would like to thank Bob Haimes, Woody Hoburg, and Berk Ozturk for their input into the technical matter, along with the EnCAPS Technical Monitor Ryan Durscher.

\end{acknowledgements}
\bibliographystyle{spmpsci}      
\bibliography{bib_deck}   

\end{document}